\documentclass{article}

\usepackage{arxiv}

\usepackage{hyperref} 
\usepackage{url} 
\usepackage{booktabs} 
\usepackage{amsfonts} 
\usepackage{nicefrac} 
\usepackage{microtype} 
\usepackage{lipsum}

\usepackage{graphicx}
\usepackage{amsmath}

\usepackage{xcolor}

\usepackage{ulem}

 \usepackage[T2A]{fontenc}
\usepackage[english]{babel}

\title{Processing of optical signals by "surgical" methods for the Gelfand-Levitan-Marchenko equation}

\author{
 Sergey Medvedev$^{1,2,*}$, Irina Vaseva$^{1,2}$, Mikhail Fedoruk$^{2,1}$\\
$^{1}$ Federal Research Center for Information and Computational Technologies,\\ Novosibirsk
630090, Russia,\\
$^{2}$ Novosibirsk State University, Novosibirsk 630090, Russia,\\
* Corresponding author: medvedev@ict.nsc.ru
}

\begin{document}
\maketitle

\begin{abstract}
We propose a new method for solving the Gelfand-Levitan-Marchenko equation (GLME) based on the 
{block version of} the Toeplitz Inner-Bordering (TIB) with an arbitrary point to start the calculation. This makes it possible to find solutions of the GLME at an arbitrary point with a cutoff of the matrix coefficient, which 
{allows} to avoid the occurrence of numerical instability and to perform calculations for soliton solutions spaced apart in the time domain. 
Using 
{an} example of two solitons, we demonstrate our method and its range of applicability. An example of eight solitons shows how the method can be applied to a more complex signal configuration.
\end{abstract}


\section{Introduction}

The Zakharov-Shabat (ZS) system arises in many physical applications \cite{lamb1980elements,akulin2005coherent}. For optics, the main applications in which the ZS system is used are fiber Bragg gratings \cite{kashyap2009fiber} and nonlinear Fourier transform (NFT) \cite{turitsyn2017nonlinear}, which is a promising basis for data transmission in fiber-optic communication lines. For the ZS system, it is necessary to solve both the direct spectral problem and the inverse one, for which it is necessary to reconstruct the potential from the spectral data. In this work, we concentrate on the inverse problem and the main application for us is NFT. Comparison of the efficiency of numerical methods for Bragg lattices is given in \cite{buryak2009comparison}. There are also several basic methods for the numerical solution of the inverse problem: the solution of the Riemann-Hilbert problem \cite{wahls2016fast,kamalian2020full}, Darboux
{'s} method for purely discrete spectrum \cite{aref2016control} and a combination of continuous and discrete spectra \cite{aref2018modulation} and the TIB method for the GLME \cite{belai2007efficient}. 
We will consider exclusively the numerical solution of the GLME, because this method, in comparison with the method based on factorization for the Riemann-Hilbert problem, is more accurate \cite{wahls2016fast}.  
Darboux's method for a combination of continuous and discrete spectrum \cite{aref2018modulation} requires solving the GLME for the continuous part of the spectral data. The iterative procedure for adding discrete values leads to an increase in computational errors.
Separately, we highlight the methods that consider the inverse NFT as a twin of the forward NFT, similar to the usual Fourier transform. Therefore, to implement inverse NFT, time-reversed forward NFT algorithms can be used (see, for example, \cite {wahls2016fast,Yousefi2020}).
Among other methods for the direct solution of the GLME, we also note an algorithm based on the transition to a system of partial differential equations  \cite{Xiao2002},  integral layer-by-layer recovery method  \cite{Rosenthal2003} and the algorithm for parameterizing kernels by polynomials presented in \cite{Ahmad1998}.  
A detailed overview of methods for solving the direct and inverse NFT problem can be found, for example, in \cite{turitsyn2017nonlinear}.

The method we propose can be used to solve the GLME for a purely continuous spectrum and we also show how to use our method to find solutions for the discrete spectrum and what computational problems arise in this case.  
The words "surgical methods"\, are used in the title of the work, because we start calculations at a fixed point, as if we were making an incision, then we cut off the "extra"\, matrix element and then move left or right, depending on the version of the GLME. All this resembles surgical manipulations.

\section{Gelfand-Levitan-Marchenko equations}

The inverse problem of recovering the potential $q(t)$ for the 
{ZS} system
\begin{equation}
\psi_{1t}+i\zeta \psi_2=q(t)\psi_2,\quad \psi_{2t}-i\zeta\psi_2=\mp q^*(t)\psi_1,   
\end{equation}
where $t\in\mathbb{R}$, $q(t),\psi_1(t),\psi_2(t)\in\mathbb{C}$, $\zeta=\xi+i\eta\in\mathbb{C}^+$,
by left or right spectral data \cite{turitsyn2017nonlinear}
\begin{equation}
\Sigma_l=\left\{l(\xi),\left[\zeta_n,l_n\right]_{n=1}^N\right\},\quad 
\Sigma_r=\left\{r(\xi),\left[\zeta_n,r_n\right]_{n=1}^N\right\},
\end{equation}
where $N$ is the number of the solitons in the signal $q(t)$,
is reduced to solving the left system of integral equations
\begin{equation}\label{A1}
\begin{array}{l}
A_1^*(t,s)+\int\limits_{-\infty}^t\,A_2(t,t')\,\Omega_l(t'+s)\,dt'=0,\quad t\geq s,\\
\mp A_2^*(t,s)+\Omega_l(t+s)+\int\limits_{-\infty}^t\,A_1(t,t')\,\Omega_l(t'+s)\,dt'=0,
\end{array}
\end{equation}
or the right system
\begin{equation}\label{B1}
\begin{array}{l}
 B^*_2(t,s)+\int\limits^{\infty}_t B_1(t,t')\Omega_r(t'+s)dt'=0,\quad t\leq s,\\
\mp B_1^*(t,s)+\Omega_r(t+s)+\int\limits^{\infty}_t B_2(t,t')\Omega_r(t'+r)\,dt'=0,
\end{array}
\end{equation}
where the kernels are defined for all real $t$ by
\begin{equation}\label{Omega_l}
\Omega_l(z)=\frac{1}{2\pi}\int\limits_{-\infty}^\infty l(\xi)e^{-i\xi z}\,d\xi -i\sum\limits_{n=1}^N\,l_{n}e^{-i \zeta_n z},
\end{equation}
\begin{equation}\label{Omega_r}
\Omega_r(z)=\frac{1}{2\pi}\int\limits_{-\infty}^\infty r(\xi)e^{i\xi z}\,d\xi -i\sum\limits_{n=1}^N\,r_{n}e^{i \zeta_n z}.
\end{equation}

After solving the system (\ref{A1}) or (\ref{B1}), the potential $q(t)$ is restored by the formulas
\begin{equation}\label{q=2A2}
q(t)=-2A_2^*(t,t)=2B_1(t,t).
\end{equation}

In 
{the} equations (\ref{A1}) and (\ref{B1}), the upper sign (minus) corresponds to the situation with the presence of a discrete spectrum, the lower sign (plus) - the absence of a discrete spectrum. 


The further presentation is given only for the left GLME, since for the right GLME the transformation is similar. 
Further we omit the subscript $l$ in the notation of the kernel $\Omega_l$.

Let's make the change of variables: $k=t-s$, $n=t$,
\begin{equation}\label{knts}
X_1(k,n)=A_1(t,s),\quad X_2(k,n)=\mp A_2^*(t,s).
\end{equation}
Taking into account the condition $t\geq s$, we obtain that the functions $X_{1,2}(k,t)$ 
are defined on the half-space $k\geq 0$.
Then the equations (\ref{A1}) take the form
\begin{equation}\label{X1}
\begin{array}{l}
X_1(k,t)\mp\int\limits_0^{\infty}\,\Omega^*(2t-k-p)X_2(p,t)\,dp=0,\\
X_2(k,t)+\Omega(2t-k)+\int\limits_0^{\infty}\,\Omega(2t-k-p)X_1(p,t)\,dp=0.
\end{array}
\end{equation}
For fixed $t$ the integral operators are Hankel, since they depend only on the sum $k+p$.

Let's make the first approximation. It consists in replacing the infinite domain of integration with a finite one $ [0, P] $ for a sufficiently large value $P$. We obtain integral equations with finite limits ($0\leq k\leq P$)
\begin{equation}\label{X1P}
\begin{array}{l}
X_1(k,t)\mp\int\limits_0^{P}\,\Omega^*(2t-k-p)X_2(p,t)\,dp=0,\\
X_2(k,t)+\Omega(2t-k)+\int\limits_0^{P}\,\Omega(2t-k-p)X_1(p,t)\,dp=0.
\end{array}
\end{equation}

To obtain the Toeplitz integral operators, we make the replacement
$X_2(p,t)=Y_2(P-p,t)$.
This allows us to transform the equations to the following form ($0\leq k,m\leq P$)
\begin{equation}\label{Y1P}
\begin{array}{l}
X_1(k,t)\mp\int\limits_0^{P}\,\Omega^*(2t-P-k+p)Y_2(p,t)\,dp=0,\\
Y_2(m,t)+\Omega(2t-P+m)+\int\limits_0^{P}\,\Omega(2t-P+m-p)X_1(p,t)\,dp=0.
\end{array}
\end{equation}

To preserve the Toeplitz structure in the finite-difference approximation of the equations (\ref{Y1P}), we use a grid with a constant step size $h$. Since after transformations the potential is determined as $q(t)=\pm 2Y_2(P,t)$, the point $P$ must be included in the computational grid. To approximate both integrals, we used the simplest 
{right Riemann sum
}. 
It is this choice that gives the second order of approximation in terms of the grid step size $h$ without using the trapezoidal method. This result was confirmed numerically. 
{If $Y_2(p,t)$  is approximated by the right Riemann sum, while $X_1(p,t)$ is approximated by the left Riemann sum, we only have the first order of approximation.} 

As a result, we obtain a system of linear equations with Toeplitz matrices
\begin{equation}\label{eq1sys}
\begin{array}{l}
\vec{X}_1 \mp h T^*\vec{Y_2}=0,\\
h T\vec{X}_1+\vec{Y_2}=\vec{F},
\end{array}
\end{equation}
where vectors have dimension $M$
$$\vec{X_1}=\left(X_1(h,t),X_1(2h,t),\cdots, X_1(Mh,t)\right)^T,$$ $$\vec{Y_2}=\left(Y_2(h,t),Y_2(2h,t),\cdots, Y_2(Mh,t)\right)^T,$$ $$\vec{F}=\left(\Omega(2t-P+h),\Omega(2t-P+2h),\cdots,\Omega(2t-P+Mh)\right)^T,$$
and the elements of the Toeplitz matrix $T$ are expressed in terms of the kernel $\Omega$
$$T_{ij}=\Omega(2t-P+(i-j)h). $$
The $T$ matrix is complete, in contrast to the \cite{belai2007efficient}, in which only the lower half of the matrix is used. 
{This} makes it possible to obtain the Toeplitz form of the system matrix (\ref{eq1sys}). In our case, we have a block matrix
\begin{equation}
T_b=\left[\begin{array}{cc}E&\mp hT^*\\hT&E\end{array}\right]   
\end{equation}
with Toeplitz blocks. To solve the system (\ref{eq1sys}) numerically, we converted the system to the block-Toeplitz form with $2\times 2$ blocks and used the block version of Levinson's algorithm. This algorithm allows us to solve the system at a fixed point $t$ and, similarly to the procedure described in \cite{belai2007efficient}, allows us to find solutions of the GLME at the nearest point $t+h/2$. Since in our approach we can start calculations at an arbitrary point $t$ and then find solutions to the right of it with minimal computational costs, we named our method Generalized Toeplitz Inner-Bordering (GTIB) method. 
{Similar} calculations for the right GLME allow finding solutions to the left of the starting point $t$.

\section{Numerical experiments}

We recover a potential $q(t)$ defined on a uniform grid of the interval of length $L$ with a step size $\tau = L/M$. In general, the system (\ref{eq1sys}) is solved on the interval $[0, P]$ divided into $M$ subintervals of
width $h = 2\tau$.

Figure 1 presents the numerical results for a single soliton potential $q(t)$ defined by the spectral data as following
\begin{equation}\label{exact1sol}
q(t) = 2\eta \mbox{sech}(2\eta t - \delta)\exp[-i(2\eta t +\theta)],
\end{equation}
where $2\eta$ is a soliton amplitude, $-2\xi$ is its frequency, $\theta$ is its phase, and $\delta$ specifies the soliton position  \cite{turitsyn2017nonlinear}.
\begin{figure}[ht]
\centering
\includegraphics[trim=0.5cm 1.5cm 3cm 0.5cm,clip=true, width=0.7\textwidth, draft=false]{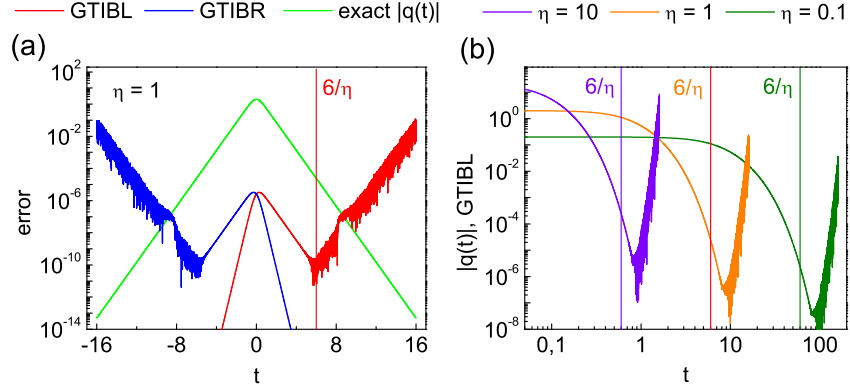}
\caption{An example of recovering a single soliton.}
\label{fig:1}
\end{figure}

Figure 1 (a) shows a single soliton $q(t)$ (\ref{exact1sol}) with $\eta = 1, \xi = 0.5, \theta = 0.8, \delta = 0$ and the errors of its recovery by GTIBL and GTIBR. L and R mean the solution of the left and right GLME correspondingly. Figure 1 (b) demonstrates numerical solutions of the left GLME for different soliton amplitudes ($\eta = 0.1, 1, 10$). We can see, that the solution stability zone depends on the amplitude of the soliton, and the boundary of the zone can be set at a distance of about $6/\eta$ from the soliton center.

Let us consider the simplest case of a two-soliton potential with parameters $\eta_1 = 1, \xi_1 = 0.5, \theta_1 = 0.1, \delta_1 = -\delta$ and $\eta_2 = 1.75, \xi_2 = -1.4, \theta_2 = 0.8, \delta_2 = \delta$. With $\delta = 8, 16, 32$ we have three two-soliton signals (see Fig. 2 (a)). The analytical formula for such a signal can be found in \cite{taha1984analytical}. It can also be determined numerically by the Darboux's method \cite{aref2016control}.  
Figure 2 (b, c, d) presents in red the errors of recovering these potentials using a combination of GTIBL and GTIBR, namely: GTIBL for $t < 0$ and GTIBR for $t >0$. Let us denote such an algorithm as "GTIB no cuts". We can see that if the solitons are too far apart (Fig. 2 (d)), the recovery error of this method is too high. 
\begin{figure}[htb]
\centering
\includegraphics[trim=0.7cm 2.8cm 3cm 0.5cm,clip=true, width=0.7\textwidth, draft=false]{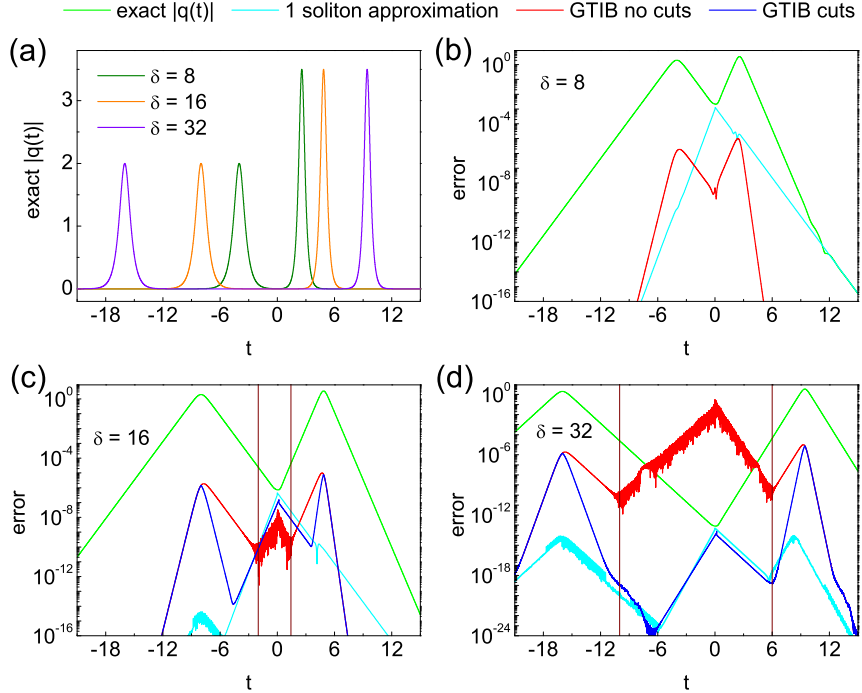}
\caption{An example of recovering two solitons located at different distances from each other.}
\label{fig:2}
\end{figure}

In order to solve this problem, we propose to cut off the matrix elements that have reached the boundary of the soliton stability zone.
To estimate how much this procedure affects the result, we can compare the exact two-soliton potential \cite{taha1984analytical} with a potential, which is specified as a combination of two single solitons (\ref{exact1sol})  (i.e. for $t < 0$ the potential is specified as the first soliton, and for $t > 0$ it is specified as the second one). The error of this one-soliton approximation of a two-soliton potential is shown in Fig. 2 (b, c, d) by cyan line.

Figure 2 (b) demonstrates the case of closely located solitons  ($\delta = 8$), their stability zones intersect. The solution obtained by using a combination of GTIBL and GTIBR gives a good accuracy (red line), while the error of one-soliton approximation, on the contrary, is large (cyan line). In this case, the cutting procedure does not make sense.

In Figure 2 (c, d) the stability zones of the solitons do not intersect. The boundaries of the stability zones are shown by brown vertical lines, they correspond to the distance $6/\eta$ from the soliton centers. In Fig. 2 (c, d) the blue line shows the error of using GTIB with cuts. This solution is obtained as follows:
from the left boundary of the computational domain to the center of the left soliton, the potential is recovered using GTIBL; from the right boundary to the center of the right soliton, the potential is recovered using GTIBR. From the midpoint $(t = 0)$ to the center of the left soliton, the potential is recovered using GTIBR, while the right soliton is cut off. Similarly, from $t = 0$ to the center of the right soliton, we use GTIBL, with the left soliton cut off.

Figure 2 (c) shows an intermediate case, when the solitons stability zones no longer intersect, but the error of GTIB without cuts is still quite small and in the central part it gives even a slightly better result than GTIB with cuts. If the solitons are located far enough, as in Fig. 2 (d), we definitely need to use the cutting procedure.

The proposed approach can be generalized from the case of two solitons to an infinite sequence of solitons. The known spectral data can be used to determine the centers of solitons and their stability zones. In each stability zone, it is possible to recover the potential by the combination of GTIBL and GTIBR, taking into account only stable solitons, and cutting off the rest. The problem arises only if each next soliton is located close to the previous one (as in Fig. 2(b)), but the edge solitons of such a sequence are located far from each other.
\begin{figure}[htb]
\centering
\includegraphics[trim=0.7cm 2.8cm 3cm 0.5cm,clip=true, width=0.7\textwidth, draft=false]{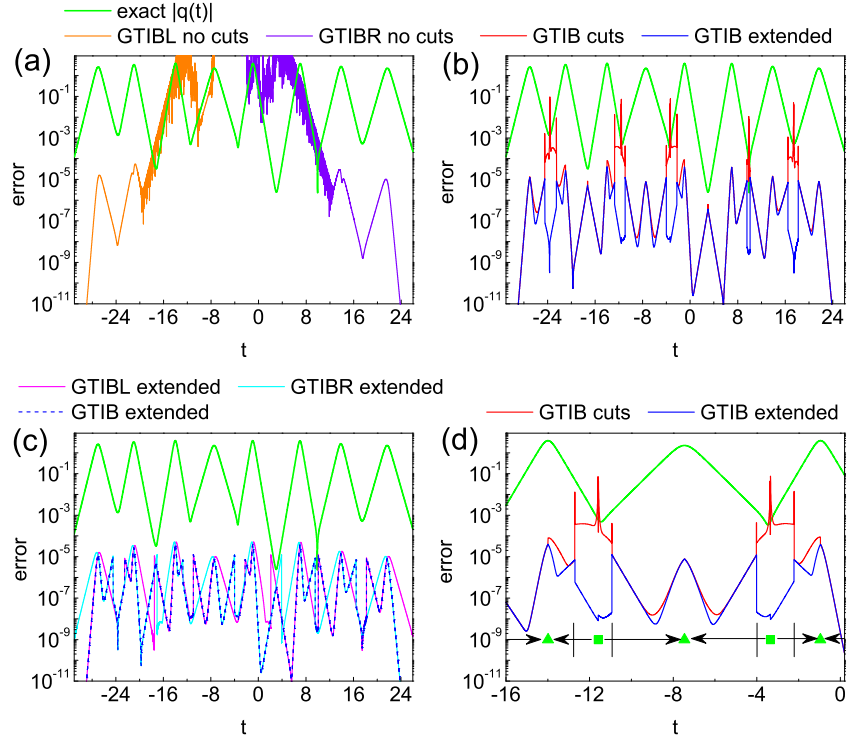}
\caption{An example of recovering eight-soliton potential.}
\label{fig:3}
\end{figure}

An example of such a situation is shown in Fig. 3 for an eight-soliton potential obtained by the Darboux's method. Fig. 3 (a) shows that a simple combination of GTIBL and GTIBR without cuts is not enough to calculate such a “long” sequence of solitons. In this case, it is impossible to avoid cutting off solitons.
Similarly to the case of a two-soliton signal, the computational domain is divided into stability zones. Calculations start at a midpoint between adjacent solitons and head left and right to their centers. 
This scheme and the locations of the cuts are shown in an enlarged fragment of the graph in Fig. 3 (d). When the cutting point is reached, the algorithm stops, cuts off the unstable soliton, and then starts again. In Fig. 3 (b, d) the red line shows an error obtained when recovering the potential in the described way.
 
A feature of the GTIB method, as well as TIB  \cite{belai2007efficient}, is a specially defined triangular domain of integration. It allows one to quickly find a solution of the GLME at a given interval, computing the next value through the previous one.
This approach also has a negative aspect, namely: the need to start from the point at which the values of the potential and matrix elements are small. Figure 3 (b, d) demonstrates that this condition is not always met, so the error of GTIB with cuts in some places turns out to be extremely large. 

We propose a modification of the method (denoted in Fig. 3 as GTIB extended), that allows us to find a solution of the GLME, starting from an arbitrary point. We expand the domain of integration in such a way that we can get more integration steps for the starting point of the algorithm. The results are shown in Fig. 3 with a blue line.

The algorithm can be simplified by using only the left (or only right) GTIB. In this case we always move in the same direction from one cut to another. The results of this approach are shown in Fig. 3 (c) in cyan and pink. It can be seen that in this case there are always zones in which the error is greater than in the case of GTIB extended, where a combination of right and left GLME is used.

We should note that if the solitons are located too close, it becomes impossible to determine their centers and to separate them from each other. If such a group is "short" (in a sense that the stability zones of its edge solitons intersect), then it can be recovered using simple combination of GTIBL and GTIBR. If such a group is "long" the potential cannot be recovered.
Nevertheless the proposed method works great for the infinite sequence of either single solitons, the stability zones of which do not intersect (see Fig. 2 (d)), or  "short" groups of close solitons (see Fig. 2(b)). Such "short" groups or single solitons can be guaranteed to be cut off from the rest of the sequence. The method still works for "long" soliton sequences as in Fig. 3, but in this case it is necessary to use the extended version of GTIB. 
\begin{figure}[h]
\centering
\includegraphics[trim=0.5cm 1.8cm 3cm 1.2cm,clip=true, width=0.6\textwidth, draft=false]{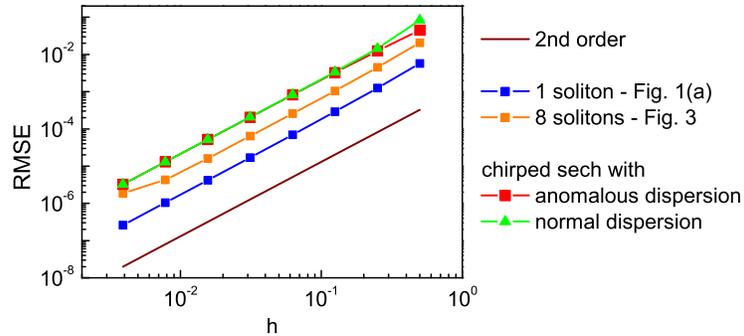}
\caption{The error of recovering different potentials with respect to the integration step size.}
\label{fig:4}
\end{figure}

Figure 4 shows the potential recovery error with respect to the integration step size $h$. A single soliton, as in Fig. 1 (a), is recovered using a combination of GTIBL and GTIBR. Eight-soliton potential, as in Fig. 3, is calculated using GTIB extended.
The potential in the form of a chirped hyperbolic secant $q(t) = A [\mbox{sech} (t)]^{1+iC}$ for $A = 5.2$, $C = 4$ is recovered using a combination of GTIBL and GTIBR. In the case of anomalous dispersion, such a signal has a continuous spectrum and five discrete eigenvalues. In the case of normal dispersion, it only has a continuous spectrum, so the matrix elements in (\ref{Omega_l}) and (\ref{Omega_r}) are bounded functions.  For all examples, the second order of approximation in $h$ is obtained.

Potential recovery errors in Fig. 1-3 and the root-mean-square error in Fig. 4 were calculated by the formulas:
\begin{equation}\label{err}
\epsilon(t)\!=\!\frac{|q(t) - q^{exact}(t)|}{\max|q^{exact}(t)|},\quad
\mbox{RMSE}[q(t)]=\sqrt{\frac{1}{M+1}\sum_{j=0}^{M}\epsilon(t_j)^2}
\end{equation}

\section{Conclusion}

A new method is proposed for the numerical solution of the GLME. It is based on the block version of the TIB  and allows to start the calculations in an arbitrary point.  The method has the second order of approximation in the integration step size. The proposed method can be used for the case, when only a continuous spectrum exists, for the case of a pure discrete spectrum, and for their combination. 
In the case of a discrete spectrum, exponentially growing matrix elements are cut off, which makes it possible to avoid the numerical instability and perform calculations for soliton solutions spaced apart in the time domain.
The method has no restrictions on the number of solitons in the signal, in contrast to the Darboux's method.

\noindent{\it Funding.} This work was supported by the Russian Science Foundation (grant No.~17-72-30006).

\noindent{\bf Disclosures.} The authors declare no conflicts of interest.
\bibliography{biblio_2sol}

\bibliographystyle{unsrt} 


\end{document}